\documentclass[12pt]{article}
\usepackage{amsmath}
\usepackage{amsfonts}
\usepackage{graphics}
\usepackage{amsthm}
\usepackage{pictex}
\usepackage{amssymb, latexsym}
\usepackage{comment}
\usepackage{rotating}
\usepackage{color}

\numberwithin{equation}{section} \numberwithin{figure}{section}
\newcommand{\HP}{Her\-mite-Pad\'e}

\newcommand{\Oo}{\mathcal{O}}

\renewcommand{\O}{\mathcal{O}}
\renewcommand{\Re}{\mathop\textrm{Re }}

\newtheorem{theorem}{Theorem}[section]
\newtheorem{corollary}[theorem]{Corollary}

\newcounter{remark}
\renewcommand{\theremark}{\arabic{section}.\arabic{remark}}

\newtheorem{definition*}{Definition}
\begin{document}
\title{Hyperelliptic uniformization \\of algebraic curves of the third order}
\author{A. I. Aptekarev*, \,\,D. N. Toulyakov\footnote {Research supported by a grant of the Russian Science Foundation (project 14-21-00025).
} \\ {\normalsize Keldysh Institute of Applied Mathematics,}
\\          {\normalsize Russian Academy of Sciences, Moscow, Russia}
        \and
        W. Van Assche\footnote{Supported by FWO research project G.0934.13, KU Leuven research grant OT/12/073 and the Belgian
Interuniversity Attraction Poles Programme P7/18.}  \\
        {\normalsize Department of Mathematics} \\
        {\normalsize Katholieke Universiteit Leuven, Belgium}}
\maketitle

\begin{center}
Dedicated to our friend Pablo Gonz\'alez Vera
\end{center}

\begin{abstract}
An algebraic function of the third order plays an important role
in the problem of asymptotics of Hermite-Pad\'e approximants
for two analytic functions with branch points. This algebraic function
appears as the Cauchy transform of the limiting measure of the asymptotic
distribution of the poles of the approximants. In many cases this algebraic
function can be determined by using the given position of the branch points
of the functions which are approximated and by the condition that its Abelian
integral has purely imaginary periods. In the present paper we obtain a
hyperelliptic uniformization of this algebraic function. In the
case when each approximated function has only two branch points,
the genus of this function can be equal to $0$, $1$ (elliptic case) or $2$
(ultra-elliptic case). We use this uniformization  to parametrize
the elliptic case. This parametrization allows us to
obtain a numerical procedure for finding this elliptic curve and as a
result we can describe the limiting measure of the distribution of the
poles of the approximants.
\end{abstract}

\vskip .25 cm


\section{Introduction} \label{sec:1}
\subsection{Definition of \HP{} approximants and  motivation of the problem}
\label{sec:1.1}

This paper is devoted to the analysis of an algebraic curve. This curve plays an
important role in the analytic theory of \HP{} approximants. We start with
their definition.

\smallskip
Let $\vec{f} = (f_1,\ldots,f_p)$ be a vector of $p$ Laurent series near infinity
\begin{equation}  \label{eq:1.1}
   f_j(z) = \sum_{k=0}^\infty \frac{f_{j,k}}{z^k},\qquad j=1,\ldots,p.
\end{equation}
The \textit{\HP{} rational approximants} (of type II)
\[  \pi_{\vec{n}} = \left( \frac{Q_{\vec{n}}^{(1)}}{P_{\vec{n}}}, \ldots,
    \frac{Q_{\vec{n}}^{(p)}}{P_{\vec{n}}} \right) \]
for the vector $\vec{f}$ and multi-index $\vec{n} = (n_1, \ldots, n_p) \in \mathbb{N}^p$
are defined by
\[  \deg P_{\vec{n}} \leq |\vec{n}| = n_1 + \cdots + n_p\,,\qquad P_{\vec{n}}\not\equiv0, \]
\begin{equation}   \label{eq:1.2}
    f_j(z)P_{\vec{n}}(z) - Q_{\vec{n}}^{(j)}(z) =:
    R_{\vec{n}}^{(j)}(z) = \O\left(\frac{1}{z^{n_j+1}}\right), \qquad z \to \infty,
\end{equation}
where the $Q_{\vec{n}}^{(j)}$ are polynomials,
for $j=1,\ldots, p$.
This definition is equivalent to a homogeneous linear system of equations
for the coefficients of the polynomial $P_{\vec{n}}$.
This system always has a solution, but the solution is not necessarily unique.
In the case of uniqueness (up to a multiplicative constant) and in
case every non-trivial solution has full degree $\vec{n}$,
the multi-index $\vec{n}$ is called \textit{normal} and the polynomial $P_{\vec{n}}$ can
be normalized to be monic
\[  P_{\vec{n}}(z) = \prod_{k=1}^{|\vec{n}|} (z-z_{k,\vec{n}}). \]
The \HP{} approximants $\pi_{\vec{n}}$ provide the best local
(at infinity) simultaneous rational approximation of the vector
$(f_1,\ldots,f_p)$ of the Laurent series (\ref{eq:1.1}) with a common denominator. The
construction (\ref{eq:1.2}) was introduced by Hermite \cite{Herm} in
connection with his proof of the transcendence of $e$. See the
papers \cite{Nutt,Apt1,WVA}  for more details.

A key problem in the study of the analytical properties (convergence, asymptotics)
of the \textit{diagonal} (i.e., $\vec{n} = (n, \ldots, n)$) \HP{} approximants is to
determine the limiting distribution of the zeros of the common denominator
$P_{\vec{n}}$, which are the poles of the \HP{} approximants, i.e., to determine the weak-$\ast$ limit of the
discrete measures
\begin{equation}  \label{eq:1.4}
   \nu_{P_{\vec{n}}} = \frac{1}{n} \sum_{k=1}^{pn} \delta(z-z_{k,\vec{n}}), \qquad n \to \infty.
\end{equation}

A class of analytic functions $\{f_j\}_{j=1}^p$ with a finite number
of branch points  plays an important role in recent investigations. We
denote the sets of the branch points of the functions $f_j$ by $A_j,\,
\, j=1, \ldots , p$. We say that
$$
   f_j \in \mathcal{A}(\overline{\mathbb{C}} \setminus A_j), \quad  \sharp A_j < \infty,
\quad   A_j \subset \mathbb{C},   \quad  j = 1,\ldots , p
$$
if the Laurent expansion (\ref{eq:1.1}) is convergent in a neighborhood of
infinity and has an analytic continuation along any path in the complex plane
$\mathbb{C}$ avoiding the sets $A_j$ of the branch points.

In accordance with Nuttall's conjecture (see \cite{Nutt}) the limit of the zero counting measures \eqref{eq:1.4} exists
$$ \nu_{P_{\vec{n}}} \stackrel{*}{\to} \lambda, $$
and the Cauchy transform of
$\lambda$ (after analytic continuation) is  an
algebraic function of order $p+1$.  We denote this function by $h$ and we denote the algebraic Riemann surface of $h$ by  $\mathfrak{R}$. The
function $h$ has  branches $h_0(z),h_1(z),\ldots, h_p(z)$ which  at infinity behave as
\begin{equation}  \label{eq:1.5}
       h_0(z) = - \frac{p}{z} + \O(z^{-2}), \qquad 
    h_j(z) = \frac{1}{z} + \O(z^{-2}), \quad 
   j=1,\ldots, p.
     \qquad z \to \infty.
\end{equation}
In particular, the conjecture states that the branch $h_0$ (fixed at infinity) is the Cauchy transform of
the limiting zero counting measure
$$ h_0(z) = \int \frac{d\lambda(x)}{x-z}\,,
\qquad z \in  \overline{\mathbb{C}} \setminus \mbox{Supp}(\lambda). $$ 
Moreover, the conjecture describes the strong asymptotics of the \HP{}  approximants \eqref{eq:1.2} using Riemann-Hilbert boundary problems on the Riemann surface $\mathfrak{R}$.
However, a still unsolved problem is how to find
(for a general setting) the main ingredient of the Nuttall conjecture, i.e., the
Riemann surface $\mathfrak{R}$.

In \cite{Apt2} a potential theory approach  to define a corresponding Riemann
surface $\mathfrak{R}$ for the case of two functions with a finite number of
branch points was suggested. The potential theory approach
was proposed for the first time in \cite{GoRa} for real analytic
functions with real sets of branch points $A_j$ and nonintersecting intervals
$\Delta_j$, which are the convex hulls of~$A_j$
$$ A_j \subset \Delta_j \subset \mathbb{R},\qquad
 \Delta_j\cap\Delta_k=\emptyset, \quad 1 \leq j \neq k \leq p.  $$
An approach based on the notion of Nuttall's
condenser was proposed recently by Rakhmanov and Suetin in \cite{RaSu}. 

The investigation of the asymptotic behavior of the diagonal \HP{}
approximants  for two functions
\begin{equation} \label{eq:1.3}
   f_j \in \mathcal{A}(\overline{\mathbb{C}} \setminus A_j),
   \quad  A_j=\{a_j,b_j\},\,\,j=1,2
\end{equation}
was started in \cite{AptKuVA}. A typical example consists of the functions $f_j(z) =
\log((z-a_j)/(z-b_j)) , \,\,  j = 1,2$. Only the case where the dispositions of the branch points of $f_1$ and $f_2$ lead to a corresponding Riemann surface $\mathfrak{R}$ of genus zero was studied in \cite{AptKuVA}. 
 A complete characterization
of algebraic curves of the third order, and of genus zero,
which have the fixed projections $A_j=\{a_j,b_j\},\,\,j=1,2$ of the
branch points, was carried out in \cite{AptLyTu}, \cite{AptKaLyTu}, (see also \cite{Lo-Ya}).

In this paper, as in \cite{AptKuVA}, we consider the case of two pairs of branch points \eqref{eq:1.3} and, in order to describe the 
asymptotics of the  \HP{} approximants, we investigate a numerical procedure for finding the equation of the algebraic Riemann surface  $\mathfrak{R}$ of genus higher than zero.

\medskip

Throughout the paper we use the following notation, conventions and definitions. 
We denote the branches of a multi-valued function $h$  by $h_0(z),h_1(z),\ldots, h_p(z)$ and the variable $z$ always belongs to subsets of the extended complex plane $z\in\bar{\mathbb{C}}$. 
An algebraic function $h$ satisfies an irreducible algebraic equation of order $p+1$ with polynomial (in $z$) coefficients. For each $z\in\bar{\mathbb{C}}$ which is not a branch point of $h$, the solutions of this equation define different branches of $h$. The Riemann surface  $\mathfrak{R}$ of the algebraic function $h$ is the ($p+1$)-sheeted covering of the extended complex plane with a finite number of branch points on which $h$ becomes a single-valued function $h(\dot{z})$. Here the variable $\dot{z}$ always belongs to $\mathfrak{R}$. We denote by $\pi$ the natural projection  $\mathfrak{R}$ on $\bar{\mathbb{C}}$ (surjection) and $\pi(\dot{z})=z$. A defragmentation of $\mathfrak{R}$ (by cuts joining the branch points) into a collection $\{\mathfrak{R}_j\}_{j=0}^p $ of  open, disjunctive sets such that $  \mathfrak{R} := \bigsqcup_{j=0}^p \overline{\mathfrak{R}}_j $ and $  \pi(\overline{\mathfrak{R}}_j)= \bar{\mathbb{C}},\,\,\, j=0,\,\ldots\,, p $ is called a sheet structure. For $\mathfrak{R}$  with assigned sheet structure $\{\mathfrak{R}_j\}_{j=0}^p $  we denote by $ z^{(j)}$ the points   $ z^{(j)} \in \mathfrak{R}_j$,  $ \pi(z^{(j)})=z\in \bar{\mathbb{C}}$ and for a single-valued
function $g$ on $\mathfrak{R}$ we can select the global branches by setting $ g_j(z):=g(z^{(j)}),\,\, j=0,\,\ldots p$.

\medskip

From now on, we only consider the (already quite difficult) case $p = 2$.

\subsection{The function $h$ and the Abelian integral of Nuttall }
\label{sec:1.2}

Now we give a formal definition of the function $h$
which will be the main subject for our analysis in this paper.
Given four points $(a_j, b_j)_{j=1,2}$, denote by $\Pi_4$ the monic polynomial with roots at these points
\begin{equation} \label{eq:1.6}
\Pi_4(x)\,:=\,\prod_{j=1}^2 (x-a_j)(x-b_j)\,=:\,x^4-s_1 x^3+s_2x^2-s_3x+s_4 .
\end{equation}
We are looking for an algebraic function $h$ of the third order,
which satisfies several requirements:
\begin{enumerate}
\item It has poles only at the branch points of $\mathfrak{R}$ (i.e., the algebraic Riemann surface where $h$ is single-valued)
whose projections belong to the given set~$A$
\begin{equation} \label{eq:1.7(1)}
h(\dot{z})\,<\,\infty\,,\quad z\in \bar{\mathbb{C}}\setminus A\,,
\qquad A:=A_1 \cup A_2,\,\,\,A_j=\{a_j,b_j\},\,\,j=1,2.
\end{equation}
Moreover, the order of a possible pole is less than the winding number of the branch point.
\item The  branches of $h$   at infinity behave as
\begin{equation}  \label{eq:1.5(2)}
    h_0(z) = - \frac{2}{z} + \cdots, \qquad 
   h_j(z) = \frac{1}{z} + \cdots, \quad 
   j=1, 2.
  \qquad z \to \infty.
\end{equation}
\end{enumerate}
Other requirements for $h$ will be given latter. The
first two conditions imply that this function has to satisfy the
equation
\begin{equation}\label{eq:1.7}
h^3-3\,\frac{P_2(z)}{\Pi_4(z)}\,h+2\,\frac{P_1(z)}{\Pi_4(z)}\,=\,0\;,
\end{equation}
where the monic polynomials $P_1$ and $P_2$ are such that $\deg P_j=j$, $j=1,2$. Indeed, if we consider the rational function 
$H(z):=h_0(z)+h_1(z)+h_2(z)$ on $\bar{\mathbb{C}}$, then a possible pole at the point $a\in A_1\cup A_2$ has order $(z-a)^\alpha$ 
for $\alpha>-1$, and therefore the rational function $H$ is bounded in $\bar{\mathbb{C}}$, which implies that $H$ is constant in $\bar{\mathbb{C}}$. 
The behavior of $h$ at infinity shows that the constant is $0$. The analysis at infinity of the other symmetric functions of $\{h_j\}_{j=0}^2$  
gives the form of the remaining coefficients of the equation \eqref{eq:1.7}.

To find out more of the polynomials $P_1$ and $P_2$ we note that the remainder terms in \eqref{eq:1.5(2)} have the same order as in \eqref{eq:1.5}, 
i.e.,
\begin{equation}\label{eq:1.5(3)}
   h_0(z) = - \frac{2}{z} + \O(z^{-2}), \qquad  h_j(z) = \frac{1}{z} + \Oo(z^{-2}), \quad
   j=1,2. \end{equation}
Indeed, since the branch $h_0$ has no branching at the point at infinity, the absence of
branching of $h_j$, $j=1,2$ leads to \eqref{eq:1.5(3)}. On the other hand, if
the point $\infty$ is a branch point for $h_j$, $j=1,2$, then the absence of a pole (in
the local variable) at this branch point again implies \eqref{eq:1.5(3)},
i.e., having a branch point there we can exclude from the above
expansions fractional degrees of $z$ between $1$ and $2$. Furthermore \eqref{eq:1.5(3)} implies
the following estimate of the discriminant $D(z)$ of the algebraic equation of $h$ in the neighborhood of the point at infinity
$$
D(z):= \prod_{k<j} (h_k-h_j)^2 =
\Oo(z^{-8}).
$$
The discriminant of the equation (\ref{eq:1.7}) is
\begin{equation}
D:=\frac{\widetilde{D}}{\Pi_4^3}, \quad
\widetilde{D}:=P_2^3-\Pi_4P_1^2,\label{eq:1.11}
\end{equation}
and we see that the assumptions on $h$ lead to
\begin{equation}\label{degreeD}
\mbox{deg }\widetilde{D}\le 4,
\end{equation}
which gives a linear equation for the coefficients of $P_1$ and $P_2$ and we obtain the representation
\begin{equation}
P_1(z):=z-c,\qquad  
P_2(z):=z^2-\displaystyle\frac{s_1+2c}{3}z+d,
\label{eq:1.8}
\end{equation}
with two unknown parameters $(c,d)$.
We have to impose extra conditions on $h$ to find these  unknown parameters $(c,d)$.

To get an extra condition on $h$ we recall the formal definition of the Abelian integral  which plays the key role
for the general conjectures of Nuttall  (see \cite{Nutt}) on the asymptotics
of  \HP{} approximants. For an algebraic three-sheeted  $\mathfrak{R}$ the \textit{Abelian integral of Nuttall} $G$ is defined by the following two
conditions.
\begin{enumerate}
\item The Abelian integral $G$ is regular on the whole  $\mathfrak{R}$, except at
$\{\infty^{(j)}\}_{j=0}^{2}$, where it has logarithmic
singularities
\begin{equation}\label{G1}
G(\dot{z})\sim
- 2\ln z, \quad \dot{z}\to\infty^{(0)}, \qquad \quad
G(\dot{z})\sim
\ln z, \quad  \dot{z}\to\infty^{(j)},\;j=1,2.
\end{equation}
\item $G$ has purely imaginary periods on $\mathfrak{R}$.
\end{enumerate}
It is known \cite{Sigel} that for any algebraic Riemann surface $\mathfrak{R}$ the Abelian integral  of Nuttall $G$
exists and the above conditions   define it uniquely
up to an additive constant $C$.  We can write
\begin{equation}\label{G99}
 G (\zeta)  = C\,+\, \int^{\zeta}_{\infty^{(0)}} H(\dot{z}) \, d \pi(\dot{z})\,,\qquad \zeta \in \mathfrak{R},
\end{equation}
where a single-valued and meromorphic function $H$ on $\mathfrak{R}$ is the derivative of $G$.
Since the integral is defined up to a constant $C$, the choice of the starting point for the integration path does not play a role and we can put it
at an arbitrary point. We choose to put it at $\infty^{(0)}$.

The second condition from above means that the function
$g:=\Re G$ is single-valued on  $\mathfrak{R}$ up to the additive constant $\Re C$ from \eqref{G99}. To fix this constant we set
$$\left(\,g_0(z)+g_1(z)+g_2(z)\,)\right|_{z=\infty}=0.$$
Thus, the second condition is equivalent to the existence of the single-valued function $g$ on $\mathfrak{R}$ 
\begin{equation}\label{G2}
g(\dot{z}):=\Re G(\dot{z}), \quad g_j(z):=g(z^{(j)}),\quad \dot{z}\in \mathfrak{R},\,\, z^{(j)}\in \mathfrak{R}_j,\quad j=0, 1, 2.
\end{equation}

We see that, due to \eqref{eq:1.5(2)} and \eqref{G1},  the Abelian
integral $G$ of Nuttall  on the Riemann surface $\mathfrak{R}$ of
the function $h$ is represented by
\begin{equation}\label{G9}
  G (\zeta)  = C\,+\, \int^{\zeta}_{\infty^{(0)}} h(\dot{z}) \, d z\,,\qquad \zeta \in \mathfrak{R},
\end{equation}
if and only if the condition \eqref{G2} is fulfilled, i.e., along any closed contour $\mathfrak{G}$ on $\mathfrak{R}$
we have
\begin{equation}
\Re \oint_{\dot{z}\in \mathfrak{G}} h(\dot{z}) \, d z= 0\;. \label{eq:1.9}
\end{equation}
Thus, we shall use the condition \eqref{eq:1.9}, i.e., the vanishing of the real
parts of the periods of the integral \eqref{G9},  as an extra condition
on $h$ in order to find the  unknown parameters $(c,d)$.

\medskip

Now we show that in our case $p=2$ and $A_j$, $j=1,2$ as in \eqref{eq:1.3} the conditions \eqref{eq:1.7(1)}, \eqref{eq:1.5(2)} and \eqref{eq:1.9} define no more than a finite number of
functions $h$, i.e.,
$\sharp\{(c,d)\}<\infty$.

Indeed, (see \eqref{eq:1.11}, \eqref{degreeD}) we have three possibilities:
\begin{equation}
\left\{\begin{array}{l}0)\ \mbox{genus } h=0 \quad\Rightarrow\quad \widetilde{D}(z)=const\,(z-z_1)^2(z-z_2)^2\;;\\
1)\ \mbox{genus } h=1 \quad\Rightarrow\quad
\widetilde{D}(z)=const\,(z-z_0)^2(z-\widetilde{z}_1)(z-\widetilde{z}_2)\;;\\
2)\ \mbox{genus } h=2\;. \label{eq:1.12}
\end{array}\right.
\end{equation}
\begin{itemize}
\item \underline{For the case $\mbox{genus } h=0$}, the condition~(\ref{eq:1.9}) is always fulfilled,
however~(\ref{eq:1.12})-0 gives two algebraic conditions (double zeros $z_1,\,z_2$) for the
determination of the two parameters $(c,d)$. In \cite{AptLyTu}, \cite{AptKaLyTu} the problem
$(a_j,b_j)_{j=1}^2\rightarrow (c,d)$ for this case was solved and it was shown that for the case
$\mbox{genus } h=0$
$$\sharp\{(c,d)\}=12.$$
\item
\underline{For the case $\mbox{genus } h=1$}, we have  in~(\ref{eq:1.12})-1 one algebraic (complex) condition
(a double zero $z_0$) and, since $\mathfrak{R}$ for this case has two cycles,
(\ref{eq:1.9}) provides two real conditions. Thus we have one complex and two real conditions
for the determination of the two complex parameters $(c,d)$.
\item
\underline{For the case $\mbox{genus } h=2$}, condition~(\ref{eq:1.9}) gives 4 real valued
relations for the two complex parameters $(c,d)$.
\end{itemize}
To make a choice of the unique $h$ from the finite set
of functions satisfying the conditions \eqref{eq:1.7(1)}, \eqref{eq:1.5(2)} and \eqref{eq:1.9} we have to analyze the 
geometrical structure of the set $\Gamma \subset \mathbb{C}$:
\begin{equation}\label{eq:1.13}
  \Gamma := \{ z \in \mathbb{C} : g_\ell(z) = g_k(z),\
        \textrm{for some } 0 \leq \ell < k \leq 2\} . \end{equation}
Thus our main goal is to start from the input data (\ref{eq:1.6})  to
get all admissible parameters $(c,d)$ and to
analyze for them the geometrical structure of the set $\Gamma $. As we already mentioned for the case of
$\mbox{genus } h=0$, this problem was completely solved in \cite{AptLyTu}, \cite{AptKaLyTu}.
In this paper we consider the cases $\mbox{genus } h >0$.

\subsection{Structure of the paper and results}
\label{sec:1.3}
We start in Section~\ref{sec:2} with some simple explicit examples of functions $h$, as in \eqref{eq:1.6}--\eqref{eq:1.8},
of genus 1 and some examples of functions $h$ of genus 2 obtained by a numerical procedure.

We have to note that the class of  functions $h$, as in \eqref{eq:1.6}, \eqref{eq:1.7}, \eqref{eq:1.8}, of genus 2 is a generic class,
in the sense that if (in the process of computation) we slightly perturb our data, then we still remain in the case of genus 2. In contrast, the case of genus 1 is a degenerated case (we have a double zero $z_0$ which bifurcates under perturbations and
moves $h$ from the case of genus 1 to the generic case of genus 2). This induces complications for the numerical procedure
of finding an appropriate $h$ of genus 1. The main goal of our paper is to find a parametrization of $c,\,d$
in  \eqref{eq:1.8}, such that the algebraic function $h$ has genus 1. This allows us to run a numeric procedure which remains in the non-generic class of genus 1.

In Section~\ref{sec:3} we find the explicit form of a conformal map of the three-sheeted Riemann surface $\mathfrak{R}$ of the algebraic function $h$, as in \eqref{eq:1.6}, \eqref{eq:1.7}, \eqref{eq:1.8}, on a two-sheeted Riemann surface $\mathfrak{H}$. It gives a hyperelliptic
(elliptic or ultra-elliptic, depending on the genus) uniformization of the algebraic curve \eqref{eq:1.7}. In the notation of
\eqref{eq:1.6}, \eqref{eq:1.7}, \eqref{eq:1.8} we have
\begin{theorem}\label{T1}
A conformal equivalence of the Riemann surfaces $\mathfrak{R}(h,z)$  and  $\mathfrak{H}(\Delta,R)$ is established by the formulas
\begin{equation}
\left\{\begin{array}{l}R = z- \displaystyle\frac{1}{h}\,,\\\\
\Delta =\displaystyle\frac{2 K_2\left(z-\frac{1}{h}\right)}{h}\,+\,\Pi_4^{\,\prime}\left(z-\frac{1}{h}\right),
\end{array}\right.\label{T1-1}
\end{equation}
where $\Pi_4^{\,\prime}(\cdot)$ is the derivative of the polynomial  \eqref{eq:1.6} with respect to its variable and $K_2(\cdot)$ is the polynomial $$K_2(x):=3x^2+2(c-s_1)\,x-3d+s_2,$$
and the inversion of \eqref{T1-1} is
\begin{equation}
\left\{\begin{array}{l} z=R+\displaystyle\frac{-\displaystyle\frac{d}{dR}\Pi_4(R)+\Delta}{2K_2(R)}, \\
h=-\,\,\displaystyle\frac{\displaystyle\frac{d}{dR}\Pi_4(R)+\Delta}{2\Pi_4(R)}.
\end{array}\right.\label{T1-2}
\end{equation}
In the new variables $(\Delta, R)$ the equation of the algebraic curve  \eqref{eq:1.7} becomes
\begin{equation}
\Delta^2\,-\, \left(\displaystyle\frac{d}{dR}\Pi_4(R)\right)^2\,+\,4\Pi_4(R)\,K_2(R)\,=\,0.
\label{ultraelliptic}
\end{equation}
\end{theorem}
\noindent  From this theorem we get an expression for the Abelian
integral of Nuttall in the variables $(R,\Delta)$.
\begin{corollary}\label{CoT1}
\begin{equation}\label{Co1}
\int h\,dz\,= \, -\, \int  \displaystyle\frac{\Delta\, dR}{2\Pi_4(R)}\,-\, \frac12 \ln \Pi_4(R)\,-\,
\ln h\,.
\end{equation}
\end{corollary}
\noindent
In Section~\ref{sec:3} we also consider an example of a conformal  mapping of $\mathfrak{R}(h,z)$ of genus 1 with branch points $\{-1, -3/8, 3/8, 1\}$ on the two sheeted Riemann surface $\mathfrak{H}(\Delta,R)$.

\medskip

Section~\ref{sec:4} is devoted to the parametrization of the algebraic curve $h$ of genus 1. As a parameter we choose the projection on the plane $R$ of the image on $\mathfrak{H}(\Delta,R)$ of the double zero $z_0$ of $\widetilde{D}$ of the discriminant \eqref{eq:1.11}.
We have
\begin{theorem}\label{T2}
If the coefficients of the equation \eqref{eq:1.7} of the algebraic function $h$ are taken in the form
\begin{equation}
\label{T2-1}
P_2(z)=\displaystyle\frac{1}{6}\,\Pi''_4(z)-\displaystyle\frac{1}{3}\,K_2(z),
\qquad
P_1(z)=\displaystyle\frac{1}{6}\,\Pi'''_4(z)-\displaystyle\frac{1}{2}\,K'_2(z),
\end{equation}
where the polynomial $K_2(\cdot)$ has the form
\begin{multline}
\label{T2-2}
K_2(x)=\frac{1}{4\Pi_4^2(R_0)}\left[12(x-R_0)^2 \Pi_4^2(R_0)+\right. \\
 \left.(x-R_0)\bigl(2\Pi_4(R_0)\Pi'_4(R_0)\Pi''_4(R_0)+
\Pi_4^{'3}(R_0)\bigr)+ \Pi_4(R_0)\Pi_4^{'2}(R_0)\right]
\end{multline}
for some $R_0\in \mathbb{C}$, then $\textup{genus } h = 1$, and conversely for any $h$ of genus 1 defined by \eqref{eq:1.7},
there exists a parameter $R_0\in \mathbb{C}$ such that the coefficients of equation \eqref{eq:1.7} have the representation \eqref{T2-2}.
\end{theorem}

Finally, in Section~\ref{sec:4} we consider an example of the application of the parametrization \eqref{T2-2}--\eqref{T2-2}. For $h$ in \eqref{eq:1.7}--\eqref{eq:1.8} we find the unknown coefficients $(c,d)$ by a numeric procedure, such that $\mbox{genus } h = 1$ and condition \eqref{eq:1.9} is fulfilled.

\section{Examples of solutions for genus $h>0$}\label{sec:2}

\subsection{The symmetric case (genus $h=1$)}\label{sec:21}
\begin{figure}[ht!]
\includegraphics[width=0.5\textwidth]{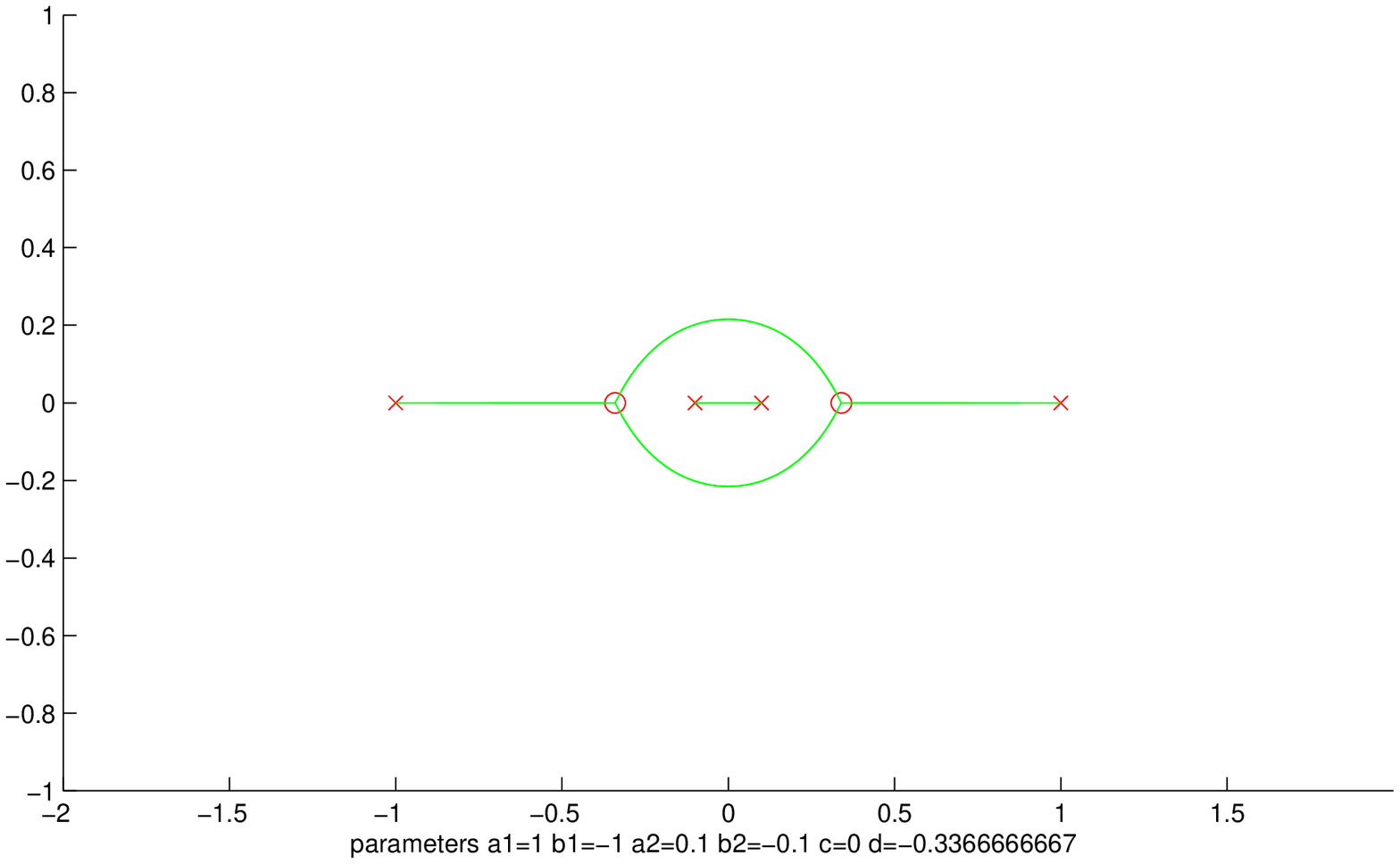}
\includegraphics[width=0.5\textwidth]{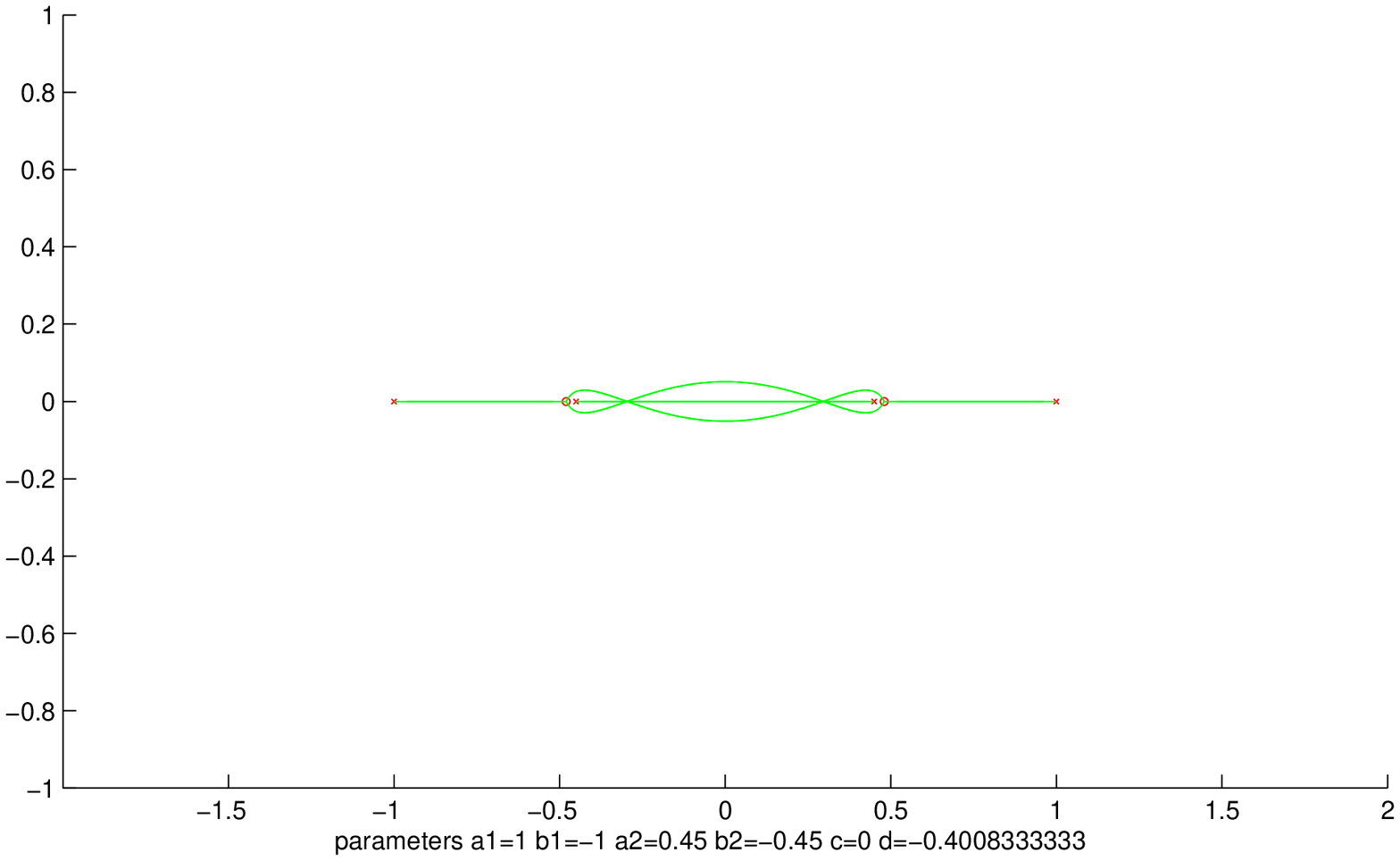}
\includegraphics[width=0.5\textwidth]{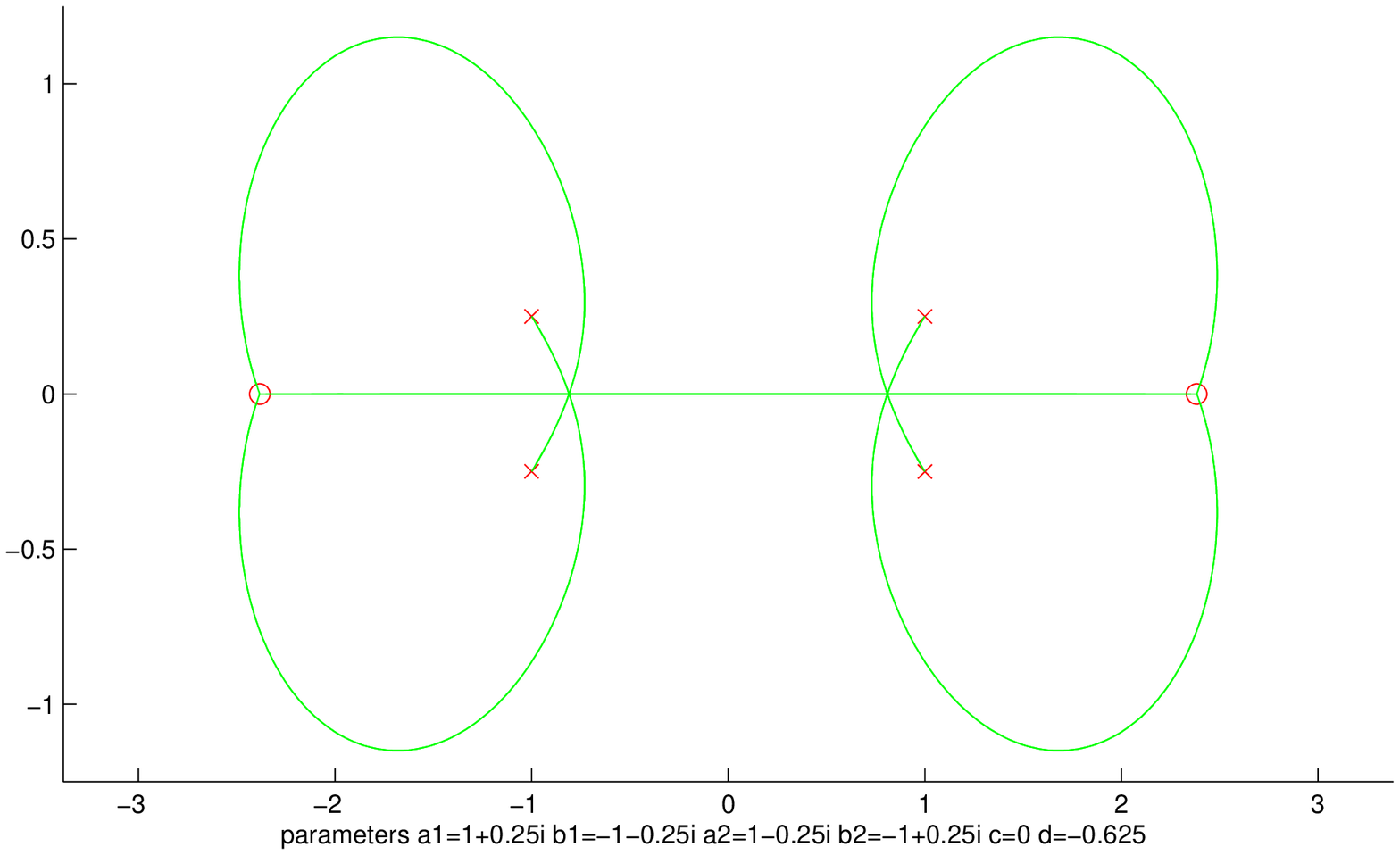}
\includegraphics[width=0.5\textwidth]{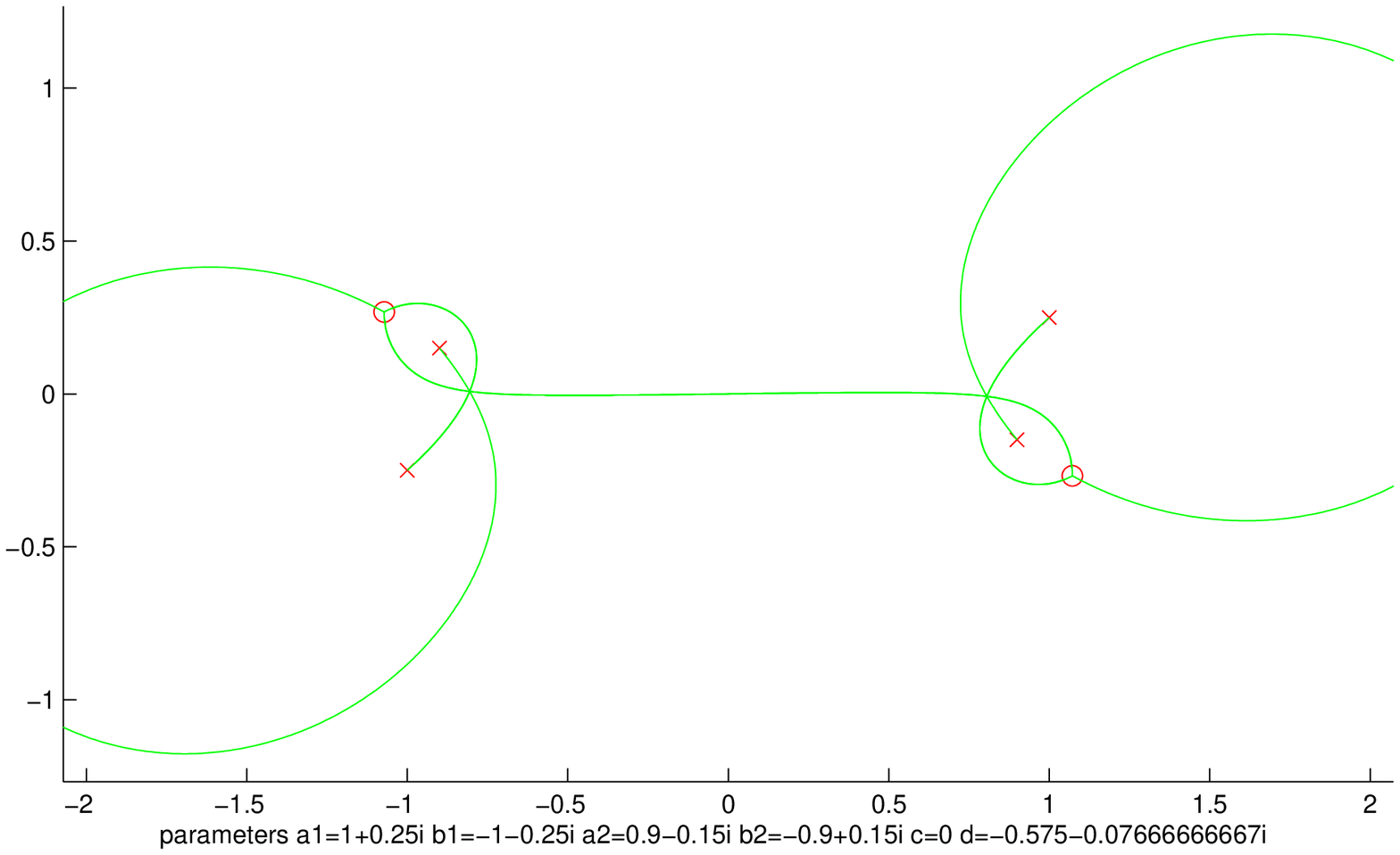}
\caption{The set $\Gamma$ for $\mbox{genus } h=1$ and symmetric input $\{a_1,b_1,a_2,b_2\}=\{a,-a,b,-b\}$.}
\label{N1}
\end{figure}

We start with some simple examples of the algebraic function $h$ of
genus 1 for which the coefficients of the equation \eqref{eq:1.7}
can be found explicitly. It is the case when the input data, i.e.,
the branch points $\{a_1, b_1, a_2, b_2\}$, have central symmetry with
respect to the origin. Indeed, for the symmetric case:
$\{a_1,b_1,a_2,b_2\}=\{a,-a,b,-b\}$ we have
$$
P_1(z):=z,\quad
P_2(z):=z^2+d,\quad
\quad \Pi_4(z)=(z^2-a^2)\,(z^2-b^2).
$$
Then, due to the symmetry, we get $c=0$ and the transcendental
condition~(\ref{eq:1.9}) is fulfilled automatically. For the
algebraic condition in (\ref{eq:1.12}-1) we have
$$
\widetilde{D}(z)=(z^2+d)^3-\Pi_4(z)\,z^2\;, \qquad \deg \widetilde{D}=4.
$$
There are two ways  to keep the symmetry: to put the double zero $z_0$
of the discriminant $\widetilde{D}$ at the point at infinity or at the
origin. In order to get the set $\Gamma$ \eqref{eq:1.13}
corresponding to the asymptotics  of the diagonal \HP{}
approximants  for \eqref{eq:1.3} we choose $z_0=\infty$. Then the
vanishing of the coefficient of $z^4$ of $\widetilde{D}$ gives
the explicit expression for the unknown parameter
\begin{equation}
\label{dsym}
d=-\frac{a^2+b^2}{3},
\end{equation}
and for the extra branch points of $h$ (the so-called \textit{soft edges})
$$
\widetilde{z}_{1,2}=\pm\frac{(a^2+b^2)^2}{3\sqrt{a^6+b^6}}.
$$
In Figure~\ref{N1} we have given  examples of the set $\Gamma$ for the
algebraic curve \eqref{eq:1.7} of genus 1 for the following
symmetric input data \eqref{eq:1.6}:
\begin{enumerate}
\item $a_1=1, b_1=-1; a_2=0.1, b_2=-0.1, d=-\frac{101}{300}$.
\item $a_1=1, b_1=-1; a_2=0.45, b_2=-0.45, d=-0.400833$.
\item $a_1=1+i\,0.25, b_1=-1+i\,0.25; a_2=1-i\,0.25, b_2=-1-i\,0.25, d=-0.625$.
\item $a_1=1+i\,0.25, b_1=-0.9+i\,0.15; a_2=0.9-i\,0.15,b_2=-1-i\,0.25, d=-0.575-i0.07666$.
\end{enumerate}
The hard edges are indicated with a cross $(\times)$ and the soft edges with a circle ($\circ$).

Regarding the relation of  $\Gamma$ and  the \HP{}
polynomials for \eqref{eq:1.3} we note the following. The first two examples for
$\Gamma$ on Figure~\ref{N1} correspond to the case when the
intervals $(a_{j},b_{j})$, $j=1,2$ are one inside the other. In case when the intervals do not intersect (the so called Angelesco case), then
$\mbox{genus } h=0$, see \cite{GoRa}, \cite{AptKuVA}. The last two examples
correspond to the case when these intervals are parallel and close
to each other. If the distance between them increases, then
again $\mbox{genus } h=0$, see \cite{AptKaLyTu}.

\subsection{``Optimization" procedure (genus $h=2$)}\label{sec:22}
B. Beckermann suggested a numerical method to find the unknown
parameters $(c,d)$ using a geometrical structure of the set
$\Gamma$. The idea of the method is simple. You start from some
approximate values $(\tilde{c},\tilde{d})$ for $(c,d)$ and from each
branch point of the corresponding  algebraic function $\tilde{h}$
you  draw the elements of the set $\Gamma$, using  the trajectories of
the quadratic differentials. Then an optimization procedure is used
to reach a coincidence of the corresponding trajectories
which started from different branch points.

This numerical method works only for functions $h$ of
genus 2, because these curves are in generic position among the curves defined
by~\eqref{eq:1.7}. The curves of lower genus are the result of algebraic
degenerations and in order to find a numerical procedure for them we need to
have their special parametrization.

\medskip

\begin{figure}[ht!]
\includegraphics[width=1.0\textwidth]{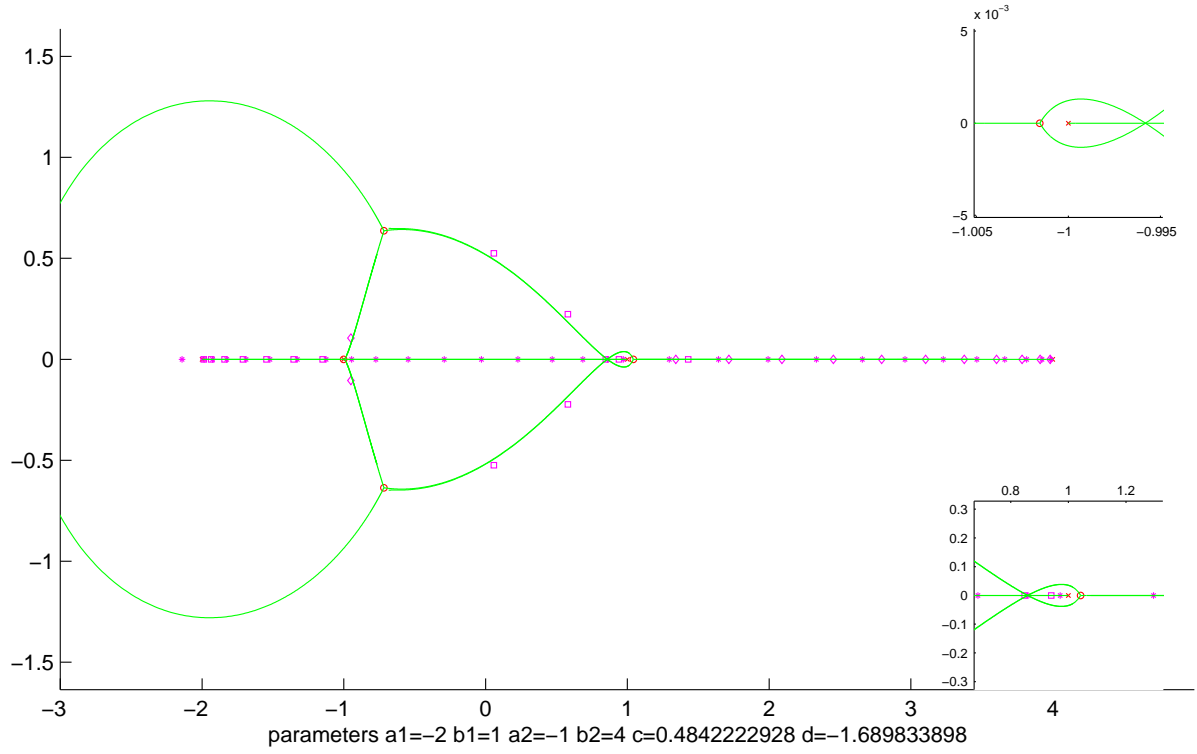}
\medskip

\includegraphics[width=1.0\textwidth]{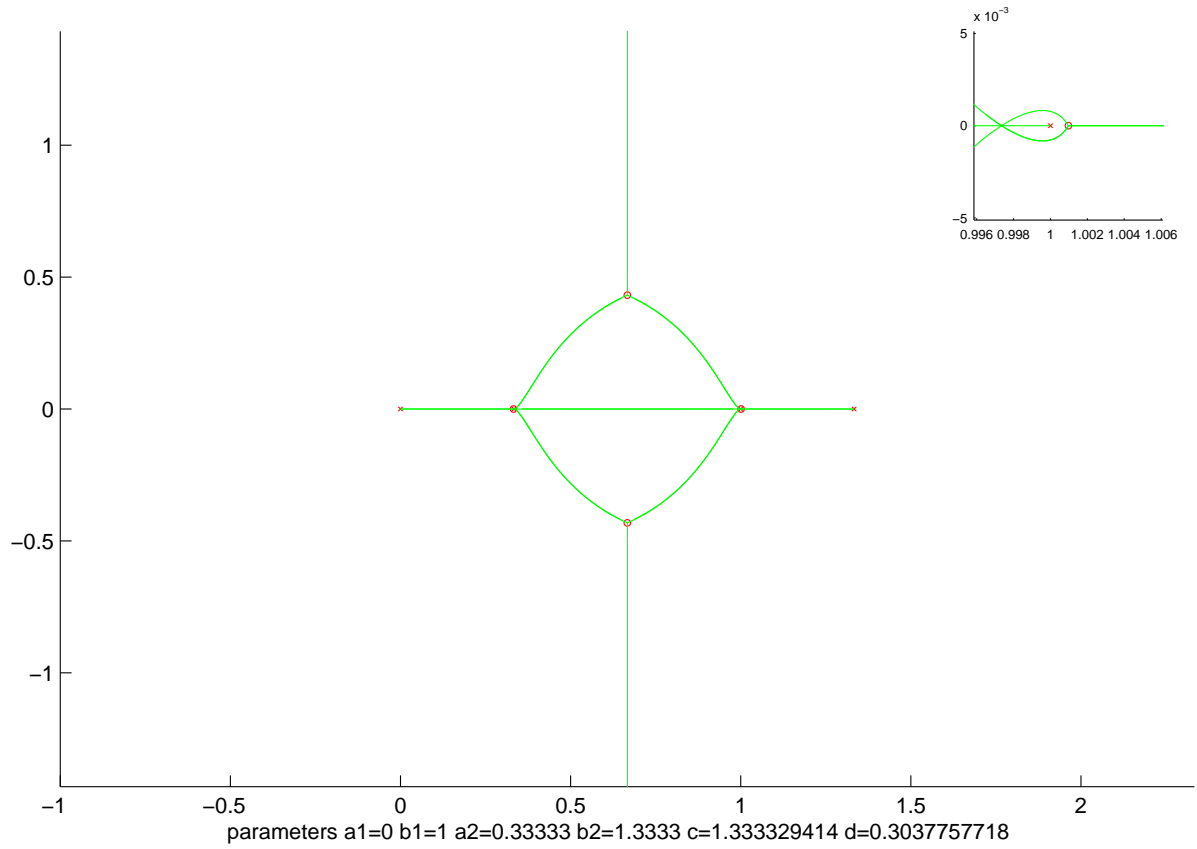}
\caption{The set $\Gamma$ for the case genus $h=2$} \label{N2}
\end{figure}

In Figure~\ref{N2} we have given  examples of the set $\Gamma$
obtained by means of this method for the following input data
\eqref{eq:1.6}:
\begin{enumerate}
\item $a_1=-2, b_1=1; a_2=-1, b_2=4; c=-0.2736665608, d=-1.689837898$.
\item $a_1=-5, b_1=1; a_2=1, b_2=5; c=0, d=-2.1647$.
\end{enumerate}
In the first example (for $\vec{n}=(20,20)$  in \eqref{eq:1.2}) we also plotted the zeros of the \HP\ polynomials of type I ({\tiny $\square$} and $\diamond$)
and type II ($\ast$).
Again the hard edges are indicated with a cross ($\times$) and the soft edges with a circle ($\circ$). The interesting
parts of the picture have been blown up to show the details.

 Observe that   H. Stahl (see \cite{St1}, \cite{St2}, \cite{AptSta}) has predicted
 $\Gamma$ (like in  Figure~\ref{N2}) for the first time for
\HP{} polynomials for functions \eqref{eq:1.3} with  overlapping
intervals~$(a_{j},b_{j})$, $j=1,2$. Another case of overlapping intervals was considered in \cite{AptLy} and \cite{Ra}.

\section{Hyperelliptic uniformization}\label{sec:3}

\textbf{Proof of Theorem~\ref{T1}.}  We are looking for  a
conformal map of the  Riemann surface of~$h$
$$
h^3-3\frac{P_2(z)}{\Pi_4(z)}\,h+2\frac{P_1(z)}{\Pi_4(z)}=0\;,
$$
on a hyperelliptic (two sheeted) Riemann surface (elliptic or
ultra-elliptic, depending on the genus of $h$). Let $Q_1$ (of exact degree one) be the divisor in the Euclidean division of $P_2$ by $P_1$ with rest $Q_0$ of degree $0$,
that is,
\begin{equation}
\label{P2}
P_2-Q_1P_1=Q_0,\quad \mbox{deg}\;Q_0=0.
\end{equation}
Then we set
$$
Q_2:=P_1Q_1^3-\Pi_4\, \Rightarrow \ \mbox{deg}\;Q_2\leqslant 2.
$$
Indeed, if we substitute $P_2$ from \eqref{P2} in \eqref{eq:1.11}, then due to \eqref{degreeD} we have
$$
\mbox{deg} \left((Q_{1} P_{1}+Q_{0})^{3} - \Pi_4 P_1^{2}\right)=
\mbox{deg}\left((Q_{2} P_{1}^{2} + 3Q_{1}^{2} P_{1}^{2}Q_{0}
 + 3Q_{1} P_{1}Q_{0}^{2} + Q_{0}^{3} \right)\leqslant 4.
$$
The last inequality implies $\mbox{deg}(Q_{2} P_{1}^{2} )\leqslant
4$ which proves $\mbox{deg}Q_{2} \leqslant 2$.

Next, we substitute the representation of $\Pi_4$ and $P_2$ by means of $Q_2,Q_1,Q_0$
and $P_1$ in the algebraic equation~(\ref{eq:1.7}) for $h$
$$
P_1Q_1^3h^3-Q_2h^3-3(P_1Q_1+Q_0)h+2P_1=0,
$$
and using the factorization $y^3+3y+2=(y-1)^2(y+2)$ for $y:=Q_1h$, we find
$$
P_1(Q_1h-1)^2(Q_1h+2)-Q_2h^3-3Q_0h=0.
$$
Now we introduce a new variable $w$ instead of $h$
$$
w=\frac{1}{h}\,,
$$
and we have
$$
P_1(z)\,(Q_1(z)-w)^2(Q_1(z)+2w)-Q_2(z)-3Q_0w^2=0.
$$
Then a new variable $R$ instead of $z$
$$
R=z-w
$$
transforms the third order equation into a quadratic equation
\begin{equation} \label{quadrEQ}
(P_1(R)+w)\,Q_1^2(R)\,(3w+Q_1(R))=Q_2(R+w)+3Q_0w^2.
\end{equation}
We note that  this magic  drop in the degree of $w$ in  \eqref{quadrEQ} is in 
correspondence with the general fact \cite{Sigel} that any compact Riemann
surface of genus less than or equal to 2 is conformally equivalent to
a two sheeted hyper-elliptic Riemann surface, i.e., ultra-elliptic, elliptic or genus zero. This fact is not generally 
valid for genus greater than 2.

Then we substitute in  \eqref{quadrEQ} the expression
$$
Q_2(R+w)= Q_2(R)+ w\,Q_2^\prime (R) + \frac{w^2}{2}\,Q_2^{\prime \prime }(R)
$$
where
$$
Q_2=P_1Q_1^3-\Pi_4\,,\quad Q_2^\prime =Q_1^3+ 3 P_1Q_1^2-\Pi_4\prime \,,\quad
Q_2^{\prime \prime} =6Q_1^2+ 6 P_1Q_1-\Pi_4^{\prime \prime}\,,
$$
and after some cancelations we arrive at
$$
\Pi_4(R)\,+\,\Pi'_4(R)\,w\,+\,\left(\frac{\Pi_4^{\prime \prime}\,-\,3P_2}{2}\right)\,w^2
\,=\,0.
$$
From this we obtain
\begin{equation}
w=\frac{-\Pi'_4(R)\pm\tilde{\Delta}(R)}{2K_2(R)}, \label{8}
\end{equation}
where
\begin{equation}
\left\{\begin{array}{l}\tilde{\Delta}^2(R):=\left[\Pi'_4(R)\right]^2-4\Pi_4(R)\,K_2(R), \\[20pt]
K_2(R):=3R^2+2(c-s_1)\,R-3d+s_2,
\end{array}\right.   \label{9}
\end{equation}
We use \eqref{8} to define the new variable
$$
\Delta =2 K_2(R) w\,+\,\Pi_4^{\,\prime}(R).
$$
Hence for each point $(z, h)$ of the Riemann surface $\mathfrak{R}$
we have a one to one correspondence with $(R, \Delta)$. If
we fix just the value of $R$, then we have a two-valued function
$\tilde{\Delta}(R)$ defined in \eqref{9}.

 For $h$ we have
\begin{equation}
h=-\frac{\Pi'_4(R)\pm\tilde{\Delta}(R)}{2\Pi_4(R)}. \label{10}
\end{equation}
Thus we have obtained a hyperelliptic uniformization for the algebraic function of the third order $h(z)$
\begin{equation}
\left\{\begin{array}{l}h=-\displaystyle\frac{\Pi'_4(R)\pm\tilde{\Delta}(R)}{2\Pi_4(R)}
\,, \\[20pt]
z=R+\displaystyle\frac{-\Pi'_4(R)\pm\tilde{\Delta}(R)}{2K_2(R)}\,.
\end{array}\right.\label{11}
\end{equation}
This proves the theorem. \qed

\bigskip

The Corollary~\ref{CoT1} of  Theorem~\ref{T1} easily
follows. Using \eqref{T1-1} we have
$$
\int h\,dz =   \int h\,dR - \ln h ,
$$
and substituting \eqref{10} we arrive at \eqref{Co1}. Another
useful representation for the Abelian integral of Nuttall by means
of elliptic integrals follows from
\begin{equation}
h\,dz = \left(\displaystyle\frac{\Pi_4'(R) K'_2(R)\,+\,4K^{2}_2(R) \,-\, 2K_2(R)\Pi_4''(R)}
{2K_2(R)\,\Delta}\,-\,\displaystyle\frac{K'_2(R)}{2K_2(R)}
\right)\,dR .
\label{hdz}
\end{equation}

\begin{figure}[h!]
\unitlength 0.6pt \linethickness{0.5pt} \centering
\framebox{\begin{picture}(200,160)(20,20)
\put(25,155){$\mathfrak{R}_0$}
\put(140,90){\circle*{4}} \put(100,90){\circle*{4}}
\put(90,90){\circle{4}} \put(60,90){\circle*{4}}
\put(150,90){\circle{4}} \put(180,90){\circle*{4}}
\put(140,92){\line(-1,0){40}} \put(140,89){\line(-1,0){40}}
\put(90,92){\line(-1,0){30}} \put(90,89){\line(-1,0){30}}
\put(180,92){\line(-1,0){30}} \put(180,89){\line(-1,0){30}}
\put(180,105){$b_1$} \put(55,105){$a_1$}
\put(130,105){$b_2$} \put(100,105){$a_2$}
\end{picture}}
\framebox{\begin{picture}(200,160)(20,20)
\put(25,155){$\mathfrak{R}_1$}
\put(90,90){\circle{4}} \put(60,90){\circle*{4}}
\put(150,90){\circle{4}} \put(180,90){\circle*{4}}
\put(90,92){\line(-1,0){30}} \put(90,89){\line(-1,0){30}}
\put(180,92){\line(-1,0){30}} \put(180,89){\line(-1,0){30}}
\put(180,105){$b_1$} \put(55,105){$a_1$}
\end{picture}}
\framebox{\begin{picture}(200,160)(20,20)
\put(25,155){$\mathfrak{R}_2$}
\put(140,90){\circle*{4}} \put(100,90){\circle*{4}}
\put(140,92){\line(-1,0){40}} \put(140,89){\line(-1,0){40}}
\put(130,105){$b_2$} \put(100,105){$a_2$}
\end{picture}}
\vskip 1 cm
\framebox{\begin{picture}(200,160)(20,20)
\put(15,165){$\mathfrak{H}_1$}
\put(70,145){\circle*{4}} \put(170,145){\circle*{4}}
\put(70,130){\small{$\epsilon_{1}$}} \put(170,130){\small{$\epsilon_{2}$}}
\put(70,35){\circle*{4}} \put(170,35){\circle*{4}}
\put(70,50){\small{$\epsilon_{3}$}} \put(170,50){\small{$\epsilon_{4}$}}
\put(140,90){\circle*{4}} \put(100,90){\circle*{4}}
\put(120,90){\circle{4}} \put(110,69){\small{$(\dot{z}=\infty^{(2)})$}}
\put(5,90){\circle{4}} \put(-30,70){\small{$(\dot{z}=\infty^{(0)})$}}
\qbezier(100,90)(107,130)(120,135)
\qbezier(140,90)(133,130)(120,135) \qbezier(100,90)(107,50)(120,45)
\qbezier(140,90)(133,50)(120,45) \put(146,105){$b_2$}
\put(80,105){$a_2$}
\end{picture}}
\framebox{\begin{picture}(200,160)(20,20)
\put(15,165){$\mathfrak{H}_2$}
\put(70,145){\circle*{4}} \put(170,145){\circle*{4}}
\put(70,130){\small{$\epsilon_{1}$}} \put(170,130){\small{$\epsilon_{2}$}}
\put(70,35){\circle*{4}} \put(170,35){\circle*{4}}
\put(60,45){\small{$\epsilon_{3}$}} \put(175,45){\small{$\epsilon_{4}$}}
\put(60,90){\circle*{4}}
\put(120,90){\circle{4}} \put(95,69){\small{$(\dot{z}=\infty^{(1)})$}}
\put(180,90){\circle*{4}}
\put(190,105){$b_1$} \put(45,105){$a_1$}
\qbezier(60,90)(65,120)(90,125)
\qbezier(90,125)(100,130)(75,165)
\qbezier(75,165)(30,185)(25,90)
\qbezier(60,90)(65,60)(90,55)
\qbezier(90,55)(100,60)(75,15)
\qbezier(75,15)(30,-5)(25,90)
\qbezier(180,90)(175,120)(150,125)
\qbezier(150,125)(140,130)(165,165)
\qbezier(165,165)(210,185)(215,90)
\qbezier(180,90)(175,60)(150,55)
\qbezier(150,55)(140,50)(165,15)
\qbezier(165,15)(210,-5)(215,90)

\end{picture}}
\caption{An example of the mapping of $\mathfrak{R}$ to $\mathfrak{H}$}
\label{R->H}
\end{figure}
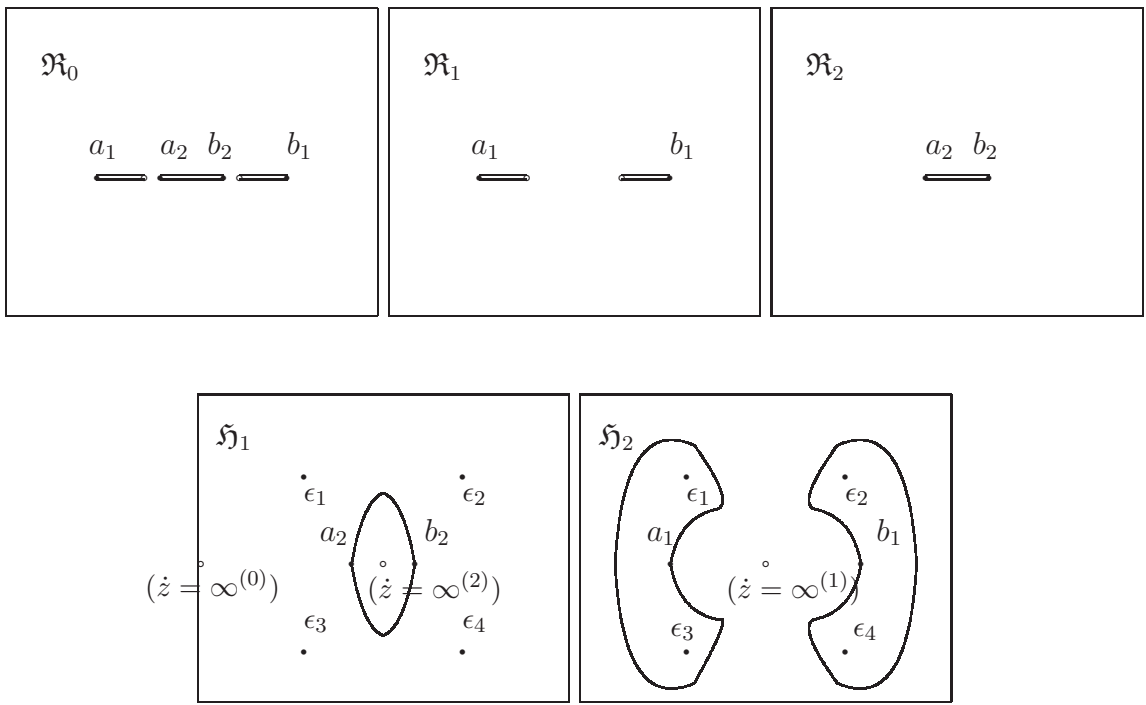

Now we consider an example of the conformal map \eqref{T1-1} of the
Riemann surface $\mathfrak{R}$ of the function $h$  of genus 1 with
branch points $\{-1, -3/8, 3/8, 1\}$ on the two-sheeted Riemann
surface $\mathfrak{H}(\Delta,R)$. Due to the symmetry we have
\eqref{dsym}, and then $d=-73/92$, $c=0$ and we can compute
$$
\tilde{\Delta}^{2}=R^{2}\,\left(4R^{4}-\frac{73}{16}R^{2}+\frac{3601}{1024}\right),\qquad
\Pi_{4}=R^{4}-\frac{73}{64}R^{2}+\frac{9}{64},\qquad K_{2}=2R^{2}.
$$
The branch points of $\mathfrak{H}$ are
$$
\tilde{\Delta}^{2}(R)=\prod_{j=1}^{4}(R-\epsilon_{j}), \qquad \epsilon_{j}= \pm
\frac{\sqrt{146 \pm 110 i \sqrt{3}}}{16}.$$
Note that the images on
$\mathfrak{H}$ of the branch points of $\mathfrak{R}$, containing the
poles of $h$,  have the same projections
 on the $R$-plane as their pre-images on the $z$-plane, i.e.,
in both cases the projections are  $A_{1}\cup A_{2}$. In
Figure~\ref{R->H} we illustrate this map, showing the images of
the sheets of $\mathfrak{R}$ and some specific values of $h(z)$.

\section{Parametrization of $h$ of genus 1.}\label{sec:4}

\textbf{Proof of Theorem~\ref{T2}.} We  have to introduce a
parameter $R_0$ and to find representations of $c(R_0)$, $d(R_0)$
for the coefficients of $P_1$, $P_2$ (see~(\ref{eq:1.8})), such
that the algebraic condition~(\ref{eq:1.12}), \eqref{eq:1.8} to obtain genus 1 for
the function $h$ is automatically fulfilled.

The dependence of the discriminant $\widetilde{D}(z)$
(see~(\ref{eq:1.11})) on $c$ and $d$ is rather complicated.
However, the appearance of the parameters $c$ and $d$ in the
discriminant $\Delta^2(R)$ (see~(\ref{9})) is linear. This easily
allows us to choose for the parameter $R_0$ the double zero of
$\Delta^2(R)$ (which characterizes the case of genus 1 for $h$)
to find representations for $c(R_0)$ and $d(R_0)$.

We denote
\begin{equation}
c_1=2(c-s_1),\quad c_2=s_2-3d.\label{12}
\end{equation}
Therefore, due to~(\ref{9}),
$$
\Delta^2(R):=\left[\Pi'_4(R)\right]^2-4\Pi_4(R)\,K_2(R)\;,\quad
K_2(R)=3R^2+c_1R+c_2\;,
$$
and the condition on $R_0$ to be a double zero of $\Delta^2(R)$ is
\begin{equation}
 \quad\left\{\begin{array}{l}\Delta^2(R_0)=0, \\[20pt]
\displaystyle\frac{d}{dR}\,\,\Delta^2(R)\Bigl|_{R=R_0}=0, \\
\end{array}\right.
\label{13}
\end{equation}
which is a linear system of equations for the determination of the parameters $c_1$ and $c_2$ (and
therefore $c,d$ as well). The solution of~(\ref{13}) is
\begin{eqnarray}
c_1 &=& \frac{1}{16\Pi_4^2(R_0)}\left(8\Pi_4(R_0)\,\Pi'_4(R_0)\,\Pi''_4(R_0)-96R_0\Pi_4^2(R_0)-4\Pi_4^{'3}(R_0)\right), \nonumber \\
c_2 &=& \frac{1}{16\Pi_4^2(R_0)}\Bigl(-8\,R_0\Pi_4(R_0)\,\Pi'_4(R_0)\,\Pi''_4(R_0)+48R_0^2\Pi_4^2(R_0) \label{14} \\
   & &  +\ 4\Pi_4^{'3}(R_0)\,R_0+4\Pi_4^{'2}(R_0)\,\Pi_4(R_0)\Bigr). \nonumber
\end{eqnarray}
In this way we have a representation for $K_2(R)$ in~(\ref{9})
\begin{multline}
 K_2(R)= \frac{1}{4\Pi_4^2(R_0)}\Bigl[12(R-R_0)^2 \Pi_4^2(R_0) \\
  +\ (R-R_0)\bigl(2\Pi_4(R_0) \Pi'_4(R_0)\,\Pi''_4(R_0) +
\Pi_4^{'3}(R_0)\bigr) + \Pi_4(R_0)\,\Pi_4^{'2}(R_0)\Bigr]. \label{15}
\end{multline}
Finally, adding the useful expressions (see~(\ref{9}) and~(\ref{eq:1.8})):
\begin{eqnarray} \label{16}
  P_2(z)&=& \frac{1}{6}\Pi''_4(z)- \frac{1}{3}K_2(z) \nonumber \\
   P_1(z)&=& \frac{1}{6} \Pi'''_4(z)-\frac{1}{2}K'_2(z),
\end{eqnarray}
we obtain the parametrization of the algebraic curve $h$
(see~(\ref{eq:1.8})) of genus 1 by the parameter $R_0$, by means of
(\ref{eq:1.7}), ~(\ref{15}) and~(\ref{16}). This proves the
theorem. \qed

\section{An example of the determination of $h$ of genus 1}\label{sec:5}
We present an  application of the $R_{0}$-parametrization in order to
determine the unknown parameters of the algebraic function
\eqref{eq:1.7}, \eqref{eq:1.8} of genus 1 for nonsymmetric input
data $\{a_1,b_1,a_2,b_2\}$.

 We consider the integral along the closed contour $\mathfrak{L}\subset \mathfrak{R}$
\begin{equation}\label{Integr}
I:=\int\limits_{\mathfrak{L}}\,  h(\dot{z}) \,d\pi(\dot{z}) .
\end{equation}
This contour (see Figure~\ref{R->H}) starts at the point $\infty^{(0)}$, goes along the pre-image of the negative part of the real axis on  
$\mathfrak{R}_0$ up to the branch point $a_1$, where it lifts to  $\mathfrak{R}_1$ and continues in the reverse direction along the pre-image of the real axis via the point $\infty^{(1)}$ up to the branch point $b_1$ where it returns to  $\mathfrak{R}_0$ and continues in the reverse direction along the pre-image of the real axis up to  the point $\infty^{(0)}$.
Then we slightly deform  $\mathfrak{L}$ so that  the new  contour $\widetilde{\mathfrak{L}}$ does not contain  the pre-images of the real axes and infinity points from the corresponding sheets of $\mathfrak{R}$. This deformation does not change the value of $I$.

We fix the parameter $R_{0}$ and substitute \eqref{15} in
\eqref{Co1}, we apply \eqref{Co1} to the integral \eqref{Integr} along $\widetilde{\mathfrak{L}}$ and we note that the periods of the 
outside integral terms are purely imaginary and therefore
  \begin{equation}\label{Integr2}
  \Re I=\Re \int\limits_{\widetilde{\mathfrak{M}}}\displaystyle\frac{\Delta(\dot{R})\, d\pi(\dot{R})}{2\Pi_4(\dot{R})}\,.
\end{equation}
Here the contour $\widetilde{\mathfrak{M}} \subset \mathfrak{H}$  is the image of  the contour $\widetilde{\mathfrak{L}} \subset \mathfrak{R}$. We introduce the notation
$$ J(R):=-\displaystyle\frac{\Delta(R)\, dR}{2\Pi_4(R)},\quad
\Delta(R):= \left(R-R_{0}\right)\,
\sqrt{\left(\displaystyle\frac{\Pi'_4(R)^2-4\Pi_4\,K_2}{(R-R_{0})^{2}}\right)},
$$
where  the continuous branch of  the square root is used for which $\sqrt{x}>0\,,\,\,\,x \rightarrow \infty$. In this notation the integral \eqref{Integr2} can be written in the form
$$
\Re I=\Re \int\limits_{\pi(\widetilde{\mathfrak{M}})} J(R) \,dR\,.
$$
  If we deform $\widetilde{\mathfrak{M}}$ to the contour $\mathfrak{M} \subset \mathfrak{H}$ which is the image of  the contour $\mathfrak{L} \subset \mathfrak{R}$ (note that $\pi(\mathfrak{M})=\mathbb{R}\subset\bar{\mathbb{C}}$ ), then the function in the integrand of \eqref{Integr2} has singularities   at the point $\dot{\infty}\in \mathfrak{M}$ and at the images of the branch points $\{a_j,\,b_j\}$ on    $ \mathfrak{M}$. 
The following behavior holds near infinity and near the branch points
$$ J(R)\sim - \frac{1}{R},\quad R\rightarrow \infty, \qquad
J(R)\sim \frac{(-1)^j}{2(R-a)},\quad R\rightarrow a\in \{a_j,\,b_j\}.
$$
Taking the regularization function
$$
\mbox{Reg} (R) := \frac{1}{2(R-a_2)}+\frac{1}{2(R-b_2)}-\frac{1}{2(R-a_1)}-\frac{1}{2(R-b_1)}
-\frac{R}{2(R^2+1)},
$$
we have 
$$\Re \int\limits_{\pi(\widetilde{\mathfrak{M}})} \mbox{Reg} (R)\,\, dR = 0.$$ 
Thus, we obtain a formula which we can use for the numeric computation of the real part of the Nuttall's Abel integral 
$$
\Re I\,=\,\Re \int\limits_{-\infty}^{\infty}(\, J(R)\,-\,\mbox{Reg} (R)\,)\, dR.
$$
Now, in accordance with \eqref{G2}, \eqref{eq:1.9} we choose $R_0$ such that $\Re I=0$.
\begin{figure}[ht!]
\includegraphics[width=1.0\textwidth]{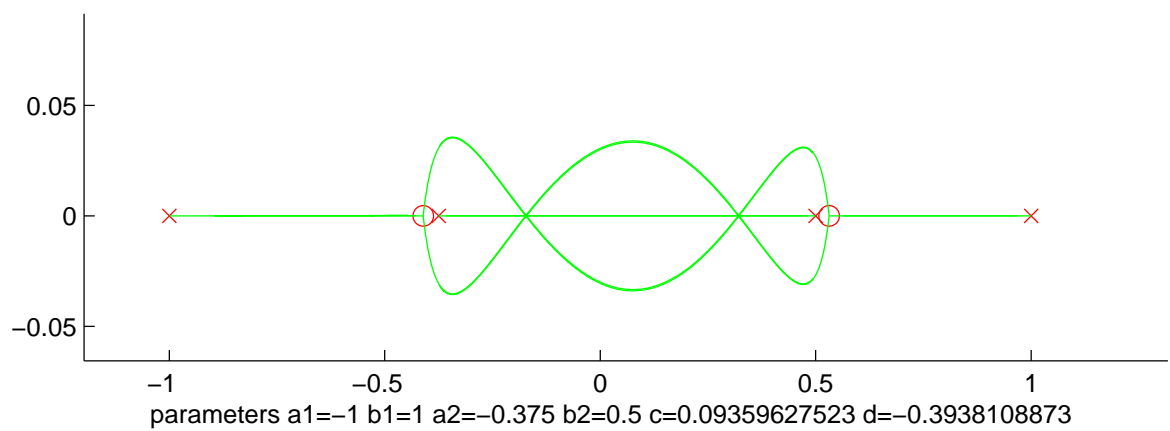}
\caption{$\Gamma$ for \mbox{genus } h\,=\,1 and symmetric input
$\{a_1,b_1,a_2,b_2\}=\{-1,-3/8,1/2,1\}$.}  \label{N3}
\end{figure}

For the input data $a_1=-1$, $b_1=1$, $a_2=-3/8$, $b_2=4/8$,  this procedure gives
$$
R_0=0.0775\quad \Rightarrow\quad \Re I=0.001,\qquad R_0=0.0774\quad \Rightarrow\quad \Re I=-0.0006.
$$
An independent verification of this procedure can be done by using the Beckermann method from Subsection~\ref{sec:22},
see Figure~\ref{N3} for $R_0=0.0775$.

\section*{Acknowledgements} 
We are very grateful to Bernd Beckermann who generously provided us with a code for computing and drawing the set $\Gamma$,
see \eqref{eq:1.13}, for the algebraic function $h$ with known coefficients
of the equations \eqref{eq:1.7} and with a code for finding these coefficients
for the generic case of genus 2.

\end{document}